\def\Aut{{\rm Aut}}
\def\Im{{\rm Im}}
\def\Ker{{\rm Ker}}
\def\F{{\bf F}}
\def\H{{\hat H}}
\def\m{{\bf m}}
\def\R{{\bf R}}
\def\Tr{{\rm Tr}}
\def\Z{{\bf Z}}
\def\Q{{\bf Q}}

\vbadness=10000
\hbadness=10000
\tolerance=10000

\proclaim Modular Moonshine III. \hfill 11 October 1994, corrected 18 Nov 1997

Richard E. Borcherds,\footnote{$^*$}{ Supported by NSF grant
DMS-9401186 and by a Royal Society research professorship.}

D.P.M.M.S.,
16 Mill Lane,
Cambridge, CB2 2SB,
England.
\bigskip

Mathematics department,
University of California at Berkeley,
CA 94720-3840,
U. S. A.
\bigskip

reb@dpmms.cam.ac.uk

www.dpmms.cam.ac.uk/\~{}reb

\bigskip

In this paper we complete the proof of Ryba's modular moonshine
conjectures [R] that was started in [B-R].  We do this by applying
Hodge theory to the cohomology of the monster Lie algebra over 
the ring of $p$-adic integers in order to
calculate the Tate cohomology groups of elements of the monster acting
on the monster vertex algebra.

\proclaim Contents. 

1. Introduction.

Notation.

2.~Representations of $\Z_p[G]$.

3.~The $\Z[1/2]$-form of the monster Lie algebra.

4.~The modular moonshine conjectures for $p\ge 13$.

5.~Calculation of some cohomology groups. 

6.~Elements of type $3B$, $5B$, $7B$, and $13B$. 

7.~Open problems and conjectures.

\proclaim
1.~Introduction.

This paper is a continuation of the earlier paper [B-R] so we will only
briefly recall results that are discussed there in more detail. In
[B-R] we constructed a modular superalgebra for each element of prime
order in the monster, and worked out the structure of this
superalgebra for some elements in the monster.  In this paper we work
out the structure of this superalgebra for the remaining elements of
prime order.

Ryba conjectured [R] that for each element of the monster of prime
order $p$ of type $pA$ there is a vertex algebra $^gV$ defined over
the finite field $\F_p$ and acted on by the centralizer $C_M(g)$ of $g$ in the
monster group $M$, with the property that the graded Brauer character
$\Tr(h|^gV)=\sum_n\Tr(h|^gV_n)q^n$ is equal to the Hauptmodul
$\Tr(gh|V)=\sum_n\Tr(g|V_n)q^n$ (where $V$ is the graded vertex
algebra acted on by the monster constructed by Frenkel, Lepowsky and
Meurman [FLM]).  In [B-R] the vertex superalgebra $^gV$ was defined for
any element $g\in M$ of odd prime order to be the sum of the Tate
cohomology groups $\H^0(g,V[1/2])\oplus
\H^1(g,V[1/2])$ for a suitable $\Z[1/2]$ form $V[1/2]$ of $V$,
and it was shown that
$\Tr(h|^gV)=\Tr(h|\H^0(g,V[1/2]))-\Tr(h|\H^1(g,V[1/2]))$ was equal to
the Hauptmodul of $gh\in M$. Hence to prove the modular moonshine
conjecture for an element $g$ of type $pA$ it is enough to prove that
$\H^1(g,V[1/2])=0$. In [B-R] this was shown by explicit calculation for
the elements of type $pA$ for $p\le 11$, using the fact that these
elements commute with an element of type $2B$. For $p\ge 13$ this
method does not work as these elements do not commute with an element
of type $2B$.  The first main theorem of this paper is theorem 4.1
which states that if $g$ is an element of prime order $p\ge 13$ not of
type $13B$ then $\H^1(g,V[1/2])=0$ (assuming a certain condition
about the monster Lie algebra, whose proof should appear later).
Hence Ryba's conjectures are
proved for all elements of $M$ of type $pA$. (Actually this is not
quite correct, because the proof for $p=2$ in [B-R] assumes an unproved
technical hypothesis.)  The proof we give fails for exactly the cases
already proved in [B-R].

We can also ask what happens for the other elements of order less than 13.
The cases of elements of types $2B$ or $3C$ are also treated in [B-R]:
for type $3C$, $\H^1(g,V[1/2])$ is again zero, and for type $2B$,
$\H^1(g,V[1/2])$ is zero in even degree and $\H^0(g,V[1/2])$ is zero
in odd degree. This leaves the cases where $g$ is of type $3B$, $5B$,
$7B$, or $13B$, when the Tate cohomology groups $\H^0(g,V[1/2])$ and
$\H^1(g,V[1/2])$ are usually both nonzero in each degree. In these
cases the structure of the cohomology groups is determined by our
second main result, theorem 6.1, which states that the unique element
$\sigma$ of order 2 in the center of $C_M(g)/O_p(C_M(g))$ acts as $1$
on $\H^0(g,V[1/2])$ and as $-1$ on $\H^1(g,V[1/2])$.  This determines
the modular characters of both cohomology groups because
$\Tr(h|\H^0(g,V[1/2])+\Tr(h|\H^1(g,V[1/2]))$ is then given by
$\Tr(h\sigma|\H^0(g,V[1/2]))-\Tr(h\sigma|\H^1(g,V[1/2]))$ which is the
Hauptmodul of the element $gh\sigma\in M$. Hence both
$\Tr(h|\H^0(g,V[1/2]))$ and $\Tr(h|\H^1(g,V[1/2]))$ can be written
explicitly as linear combinations of two Hauptmoduls for any
$p$-regular element $h\in C_M(g)$.  We get two new modular vertex
superalgebras acted on by double covers of sporadic groups when $g$
has type $3B$ or $5B$: a vertex superalgebra over $\F_3$ acted on by
$2.Suz$, and a vertex superalgebra over $\F_5$ acted on by $2.HJ$. For
elements of types $7B$ and $13B$ the vertex superalgebras we get are
acted on by the double covers $2.A_7$ and $2.A_4$ of alternating
groups.

In particular if the integral form $V$ of the monster vertex algebra
discussed in [B-R] exists and has the properties conjectured there then
all the cohomology groups $\H^i(g,V)$ are now known for all elements
$g$ of prime order in $M$.

We now discuss the proofs of the two main theorems. We first quickly
dispose of the second main theorem: if an element $g\in M$ of odd
order commutes with an element of type $2B$ then it is easy to work
out its action on the $\Z[1/2]$-form of the monster vertex algebra and
hence the Tate cohomology groups $\H^*(g,V[1/2])$ can be calculated by
brute force, which is what we do in sections 5 and 6 for the elements
of types $3B$, $5B$, $7B$, and $13B$.

The elements $g$ of prime order $p\ge 13$ do not commute with any
elements of type $2B$ (except when $g$ has type $13B$ or $23AB$) so we
cannot use the method above. Instead we adapt the proof of the
moonshine conjectures in [B]. By using a $\Z_p$-form of the monster
Lie algebra rather than a $\Q$-form we can find some complicated
relations between the coefficients of the series
$\sum_n\Tr(h|\H^0(g,V[1/2]_n))q^n$ and
$\sum_n\Tr(h|\H^1(g,V[1/2]_n))q^n$.  The reason we get more
information by using $\Z_p$-forms rather than $\Q_p$-forms is that the
group ring $\Z_p[\Z/p\Z]$ has 3 indecomposable modules which are free
over $\Z_p$, while the group ring $\Q_p[\Z/p\Z]$ only has 2
indecomposable modules which are free over $\Q_p$, and this extra
indecomposable module provides the extra information.  The relations
between the coefficients we get are similar to the relations defining
completely replicable functions, except that some of the relations
defining completely replicable functions (about ``$(2p-1)/p^2$'' of
them) are missing.  If $p$ is large ($\ge 13$) we show that these
equations have only a finite number of solutions. We then use a
computer to find all the solutions, and see that the only solutions
imply that $\H^1(g,V[1/2])=0$ for all $g$ of prime order at least 13
other than those of type $13B$. This proof could probably be made to
work for some smaller values of $p$ such as $7$ and $11$, but the
difficulty of showing that there are no extra solutions for the
equations of the coefficients increases rapidly as $p$ decreases,
because for smaller $p$ there are more equations missing.

Y. Martin [Ma] has recently found a conceptual proof that any
completely replicable function is a modular function.
Cummins and Gannon [CG] have greatly generalized this result
using different methods, and showed that the functions
are Hauptmoduls. It seems possible that their methods could be extended
to replace the computer calculations in section 4, although Cummins
has told me that this would probably require adding in some ``missing''
relations corresponding to the case $(p,mn)\ne 0$ in proposition 3.4. 

We summarize the results of this paper and of [B-R] about the vertex
superalgebras $^gV$ for all elements $g$ of prime order $p$ in the
monster. If $g$ is of type $pA$, or $pB$ for $p>13$, or $3C$, then
$^gV$ is just a vertex algebra and not a superalgebra.  If $p\ge 13$
this is proved in this paper by an argument which works for any self
dual form of the monster vertex algebra but which relies on computer
calculations and on a so far  unpublished argument about
integral forms of the monster Lie algebra.
For $p<13$ (or $p=23$) this is proved in [B-R] by
calculating explicitly for a certain $\Z[1/2]$ form of the monster
vertex algebra; this avoids computer calculations but only works for
one particular form of the monster vertex algebra.  Also, if $p=2$ the
calculation depends on an as yet unproved technical assumption about
the Dong-Mason-Montague construction of the monster vertex algebra
from an element of type $3B$ ([DM] or [M]).  If $g$ is of type $pB$
for $2\le p\le 13$ then the super part of $^gV$ does not vanish. If
$p=2$ then $^gV$ is calculated explicitly in [B-R] (again using the
unproved technical assumption) and it turns out that its ordinary part
vanishes in odd degrees and its super part vanishes in even degrees.
If $3\le p\le 13$ then $^gV$ is calculated in this paper, and we find
that there is an element in $C_M(g)$ acting as $-1$ on the super part
and as $1$ on the ordinary part of $^gV$.  In all cases we have
explicitly described the modular characters of $C_M(g)$ acting on the
ordinary or super part of any degree piece of $^gV$ in terms of the
coefficients of certain Hauptmoduls.

\proclaim
Notation.
\item{$[A]$}  is the element of $K$ represented by the module $A$.
\item{$A,B,C$} $G$-modules.
\item{$A_n$} The alternating group on $n$ symbols.
\item{$\Aut$} The automorphism group of something. 
\item{$c^+_g(n)$} The $n$'th coefficient of the Hauptmodul of $g\in M$, 
equal to $\Tr(1|^gV_n)=\dim(\H^0(g,V[1/2]_n))-\dim(\H^1(g,V[1/2]_n))$. 
\item{$c^-_g(n)$} $\dim(\H^0(g,V[1/2]_n))+\dim(\H^1(g,V[1/2]_n))$. 
\item{$c_{m,n}$} Defined in proposition 4.3. 
\item{${\bf C}$} The complex numbers. 
\item{$C_M(g)$} The centralizer of $g$ in the group $M$. 
\item{$Co_1$} Conway's largest sporadic simple group.
\item{$E$} The positive subalgebra of the monster Lie algebra.
\item{$F$} The negative subalgebra of the monster Lie algebra.
\item{$\F_p$} The finite field with $p$ elements. 
\item{$f$} A homomorphism from $K$ to $\Q$ defined in lemma 2.3.
\item {$g$} An element of $G$, usually of order $p$. 
\item {$g_i$} The element $g^i$ of $G$, used when it is necessary 
to distinguish the multiplication in the group ring from some other
multiplication.
\item{$\langle g\rangle$} The group generated by $g$. 
\item{$G$} A group, often cyclic of prime order $p$ and generated by $g$. 
\item{$h(A),h_n$} See section 5. 
\item{$H$} The Cartan subalgebra of the monster Lie algebra.
\item{$\H^i(G,A)$} A Tate cohomology group of the finite group $G$
with coefficients in the $G$-module $A$.
\item{$\H^i(g,A)$} means $\H^i(\langle g\rangle,A)$, where $\langle g\rangle$
is the cyclic group generated by $g$.
\item{$\H^*(g,A)$} The sum of the Tate cohomology groups
$\H^0(g,A)$ and $\H^1(g,A)$, considered as a super module.
\item{$HJ$} The Hall-Janko sporadic simple group (sometimes denoted $J_2$). 
\item{$I$} The indecomposable module of dimension $p-1$ over $\Z_p[G]$, 
isomorphic to the kernel of the natural map from $\Z_p[G]$ to $\Z_p$
and to the quotient  $\Z_p[G]/N_G\Z_p$. 
\item{$\Im$} The image of a map. 
\item{$K$} A ring which is a free $\Q$-module with a basis
of 3 elements $[\Z_p]=1$, $[\Z_p[G]]$,
and $[I]$.
\item{$\Ker$} The kernel of a map.
\item{$\Lambda,\hat\Lambda$} The Leech lattice and a double cover of the Leech 
lattice.
\item{$\Lambda^n(A)$} The $n$'th exterior power of $A$.
\item{$\Lambda^*(A)$} The  exterior algebra $\oplus_n\Lambda^n(A)$ of $A$.
\item{$L$} An even lattice. 
\item{$\m$} The monster Lie algebra. 
\item{$M$} The monster simple group.
\item{$N_G$} The element $\sum_{g\in G}g$ of $\Z_p[G]$. 
\item{$O_p(G)$} The largest normal $p$-subgroup of the finite group $G$. 
\item{$p$} A prime, usually the order of $g$.
\item{$\Q,\Q_p$} The rational numbers and the field of $p$-adic numbers. 
\item{$\rho$} The Weyl vector of the monster Lie algebra. 
\item{$\R $} The real numbers.
\item{$R_p$} A finite extension of the $p$-adic integers.  
\item{$R(i)$} Defined just before lemma 3.2. 
\item{$\sigma$} An automorphism of type $2B$ in the monster
or the automorphism $-1$ of the Leech lattice. 
\item{$S_n$} A symmetric group.
\item{$S^n(A)$} The $n$'th symmetric power of $A$.
\item{$S^*(A)$} The  symmetric algebra $\oplus_nS^n(A)$ of $A$.
\item{$Suz$} Suzuki's sporadic simple group.  
\item{$\Tr$} $\Tr(g|A)$ is the usual trace of $g$ on a module
$A$ if $A$ is a module over a ring of characteristic 0, and the
Brauer trace if $A$ is a module over a field of finite characteristic.
\item{$V[1/n]$} A $\Z[1/n]$-form of the monster vertex algebra. 
\item{$V_\Lambda$} The integral form of the vertex algebra of $\hat\Lambda$. 
\item{$V_n$} The degree $n$ piece of $V$. 
\item{$V^n$} An eigenspace of some group acting on $V$. 
\item{$^gV$} A modular vertex algebra or superalgebra given by $\H^*(g,V[1/n])$
for some $n$ coprime to $|g|$.
\item{$\Z$} The integers.
\item{$\Z_p$} The ring of $p$-adic integers.
\item{$\omega$} A cube root of 1 or a conformal vector.
\item{$\Omega$} The Laplace operator on $\Lambda^*(E)$.

\proclaim
2.~Representations of $\Z_p[G]$.

We give some auxiliary results about modules over $\Z_p[G]$ where $G$ is a
cyclic group generated by $g$ of prime order $p$, and in particular
calculate the exterior and symmetric algebras of all indecomposable
modules.  All modules will be free over $\Z_p$ and will be either
finitely generated or $\Z$ graded with finitely generated pieces of
each degree.  All tensor products will be taken over $\Z_p$.

Recall from [B-R section 2] that there are 3 indecomposable modules
over $\Z_p[G]$, which are $\Z_p$, the group ring $\Z_p[G]$, and the
module $I$ of $\Z_p[G]$ that is isomorphic to the kernel of the
natural map from $\Z_p[G]$ to $\Z_p$ and to the quotient
$\Z_p[G]/N_G\Z_p$ (where $N_G=1+g+g^2+\cdots+g^{p-1}$). Their Tate
cohomology groups are given by $\H^0(g,\Z_p)=\Z/p\Z$,
$\H^1(g,\Z_p)=0$, $\H^0(g,\Z_p[G])=0$, $\H^1(g,\Z_p[G])=0$,
$\H^0(g,I)=0$, and $\H^1(g,I)=\Z/p\Z$.

\proclaim Lemma 2.1. The tensor products of these modules
are given as follows:
$$\eqalign{
\Z_p\otimes X&=X \hbox{ (for any $X$)}\cr
\Z_p[G]\otimes X&=\Z_p[G]^{\dim(X)} \hbox{ (the sum of $\dim(X)$ copies of
 $\Z_p[G]$)}\cr
I\otimes I&= \Z_p[G]^{p-2}\oplus \Z_p.\cr
}$$

Proof. The case $\Z_p\otimes X$ is trivial. If $X$ has a basis
$x_1,\ldots,x_n$ then for any fixed $i$ the elements $g^j\otimes
g^j(x_i)$ ($0\le j<p$) form a basis for a submodule of $\Z_p[G]\otimes
X$ isomorphic to $\Z_p[G]$, and $\Z_p[G]\otimes X$ is the direct sum
of these submodules for all $i$, which proves the result about
$\Z_p[G]\otimes X$.  For the case $I\otimes I$ we look at the
cohomology sequence of the short exact sequence $$0\rightarrow
I\otimes I\rightarrow I\otimes \Z_p[G]\rightarrow
I\otimes\Z_p\rightarrow 0 $$ to see that $\H^0(I\otimes I)=\Z/p\Z$ and
$\H^1(I\otimes I)=0$, which implies that $I\otimes I$ must be the sum
of $\Z_p$ and some copies of $\Z_p[G]$.  The number of copies can be
worked out by looking at the dimensions of both sides.  This proves
lemma 2.1.

\proclaim Corollary 2.2. If $A$ and $B$ are $\Z_p[G]$-modules
then $\H^*(g,A\otimes B)=\H^*(g,A)\otimes\H^*(g,B)$.

Proof.  Check each of the 9 cases when $A$ and $B$ are indecomposable
modules using lemma 2.1.  This proves corollary 2.2.

We form the ring $K$ which is a free $\Q$-module with a basis of 3
elements $[\Z_p]=1$, $[\Z_p[G]]$, and $[I]$ corresponding to the
indecomposable modules of $\Z_p[G]$, whose addition and multiplication
are those induced by tensor products and direct sums of modules.  If
$A$ is a $\Z_p[G]$-module then we write $[A]$ for the element of $K$
represented by $A$.  Notice that this is not the same as the (tensor
product with $\Q$ of the) Grothendieck ring of finitely generated
projective $\Z_p[G]$ modules because we do not assume the relations
given by nonsplit short exact sequences. The Grothendieck ring
(tensored with $\Q$) is a quotient of our ring $K$ by the ideal
generated by $[\Z_p]+[I]-[\Z_p[G]]$.

\proclaim Lemma 2.3. There are exactly 3 ring homomorphisms 
$\dim$, $\Tr(g|.)$, and $f$ from $K$ to $\Q$ (which correspond to the
3 elements of the spectrum of $K$).  They are given by
\item {1.} $\dim([\Z_p])=1$, $\dim([\Z_p[G]])=p$, $\dim([I])=p-1$.
\item {2.} $\Tr(g|[\Z_p])=1$, $\Tr(g|[\Z_p[G]])=0$, $\Tr(g|[I])=-1$.
\item {3.} $f([\Z_p])=1$, $f([\Z_p[G]])=0$, $f([I])=1$.

Proof. From corollary 2.2 we know that the products in the ring $K$
are given by $[\Z_p[G]][\Z_p[G]]=p[\Z_p[G]]$,
$[\Z_p[G]][I]=(p-1)[\Z_p[G]]$, and $[I][I]= (p-2)[\Z_p[G]]+1$. Hence
if $x\in \Q$ and $y\in \Q$ are the images of $[\Z_p[G]]$ and $I$ under
some homomorphism they must satisfy $x^2=px$, $xy=(p-1)x$, and
$y^2=(p-2)x+1$.  The only solutions to these equations are
$(x,y)=(p,p-1)$, $(0,1)$, or $(0,-1)$.  This proves lemma 2.3.

The first two homomorphisms in lemma 2.3 are the homomorphisms which
factor through the Grothendieck ring and will not be very interesting to us.
We can use them in place of $f$ in the arguments in the rest of the
paper, but this only leads to relations between the coefficients of
$\sum_n\dim(V_n)q^n$ and $\sum_n\Tr(g|V_n)q^n$, which are essentially
the relations already used in [B] to determine these functions. It is
the third homomorphism $f$ which detects the integral structure on
$V[1/2]$ which will give us enough extra information to prove theorem
4.1.

In general, taking the trace of any fixed element of a group gives a
homomorphism from a representation ring of a group to some field, and
quite often (e.g., for the complex representation ring of a group)
this gives a 1:1 correspondence between homomorphisms and (certain
sorts of) conjugacy classes of the group.  Lemma 2.3 gives an example
of a homomorphism $f$ not constructed in this way.

A permutation module is a module for $G$ with a basis acted on by $G$.
We recall from [B-R] that a module $A$ is a permutation module if and
only if $\H^1(g,A)=0$, which in turn is equivalent to saying that $A$
is a sum of copies of $\Z_p$ and $\Z_p[G]$.

\proclaim Lemma 2.4. If $A$ is a permutation module over any ring
for any finite group $g$ then so are the symmetric powers $S^n(A)$.
If all elements of $G$ act on a basis of $A$ as products of cycles of
odd length then the exterior powers $\Lambda^n(A)$ are also
permutation modules.

Proof. For symmetric powers this is proved in [B-R].  For the proof for
exterior powers, let $a_1,\ldots,a_n$ be a basis for $A$ acted on by
$G$. Then $G$ acts on the set of all elements of the form $\pm
a_{i_1}\wedge\cdots\wedge a_{i_n}$, which consists of all elements of
a basis together with their negatives.  If the action of $G$ on the
basis of $A$ has the property that all elements of $G$ are products of
cycles of odd length, then there is no element of $G$ that maps any
element of the form $\pm a_{i_1}\wedge\cdots\wedge a_{i_n}$ to its
negative.  This means that we can find a $G$-invariant subset of these
elements which forms a basis for $\Lambda^n(A)$. This proves lemma
2.4.

\proclaim Lemma 2.5. The exterior  powers 
of indecomposable modules are given by
$$\eqalign{
\Lambda^0(\Z_p)&=\Z_p\cr
\Lambda^1(\Z_p)&=\Z_p\cr
\Lambda^n(\Z_p)&=0 \hbox{ if $n\ne 0,1$}\cr
\Lambda^0(\Z_p[G])&=\Z_p\cr
\Lambda^p(\Z_p[G])&=\Z_p \hbox{ if $p$ is odd, and $I$ if $p=2$}\cr
\Lambda^n(\Z_p[G])&=\Z_p[G]^{p-1\choose n} \hbox{ if $1\le n\le p-1$}\cr
\Lambda^n(\Z_p[G])&=0 \hbox{ if $n>p$}\cr
\Lambda^n(I)&=\Z_p\oplus \Z_p[G]^{({p-1\choose n}-1)/p} 
\hbox{ if $n$ is even and
	$0\le n<p$}\cr
\Lambda^n(I)&=I\oplus \Z_p[G]^{({p-1\choose n}+1-p)/p} \hbox{ if $n$ is odd and
	$0\le n<p$}\cr
\Lambda^n(I)&=0 \hbox{ if $n\ge p$}\cr
}$$

Proof.  For the exterior powers of $\Z_p[G]$ with $p$ odd we apply
lemma 2.4 to see that all exterior powers are sums of $\Z_p$ and
$\Z_p[G]$, and the exterior powers can then be identified by checking
the dimension and the trace of $g$. The cohomology groups of the
exterior powers of $I$ can then be worked out by induction using the
exact sequence $$0\rightarrow \Lambda^n(I)\rightarrow
\Lambda^n(\Z_p[G])\rightarrow
\Lambda^{n-1}(I)\rightarrow 0.$$
The cohomology groups determine how many copies of $\Z_p$ and $I$ 
occur in the decomposition into indecomposable modules, 
so we can work out the exterior powers of $I$ from this. 
The remaining cases of lemma 2.5  are trivial, so this
proves lemma 2.5. 

These exterior powers induce operators on the ring $K$ which we also
denote by $\Lambda^n$.  For example, we can use these operators to
define Adams operations $\psi^n$ on the ring $K$ by $$\sum_{n\in
\Z}(-1)^n\Lambda^n([A])q^n=\exp(-\sum_{n\in \Z}\psi^n([A])q^n/n).$$

\proclaim Corollary 2.6. If $p$ is odd then 
$$\eqalign{
\sum_{n\in\Z} (-1)^nf([\Lambda^n(\Z_p)])q^n &= 1-q\cr
\sum_{n\in\Z} (-1)^nf([\Lambda^n(\Z_p[G])])q^n &= 1-q^p\cr
\sum_{n\in\Z} (-1)^nf([\Lambda^n(I)])q^n &= (1+q^{p})/(1+q)\cr
}$$

Proof. This follows immediately from lemma 2.5 and lemma 2.3. 

\proclaim Lemma 2.7. The symmetric powers 
of indecomposable modules are given by
$$\eqalign{
S^n(\Z_p)&=\Z_p \hbox{ if $n\ge 0$}\cr
S^n(\Z_p[G])&=\Z_p[G]^{{p+n-1\choose n}/p}\hbox{ if $(p,n)=1$}\cr
S^n(\Z_p[G])&=\Z_p[G]^{({p+n-1\choose n}-1)/p}\oplus \Z_p
	\hbox{ if $p|n$}\cr
S^n(I)&=\Z_p\oplus\Z_p[G]^{({p+n-2\choose n}-1)/p}
	\hbox{ if $p|n$}\cr
S^n(I)&=I\oplus\Z_p[G]^{({p+n-2\choose n}+1-p)/p}
	\hbox{ if $n\equiv 1 \bmod p$}\cr
S^n(I)&=\Z_p[G]^{({p+n-2\choose n})/p}
	\hbox{ if $n\not\equiv 0,1\bmod p$}\cr
}$$

Proof. For symmetric powers of $\Z_p$ this is trivial and for
symmetric powers of $\Z_p[G]$ it can be proved in the same way as in
lemma 2.5. If $p=2$ then the case of $I$ is easy to do as $I$ is just
one dimensional.  For the case of $I$ when $p$ is odd we let $N_G$ be
the element $1+g+\cdots +g^{p-1}\in \Z_p[G]$ and consider the quotient
of the symmetric algebra $S^*(\Z_p[G])$ by the ideal $(N_G)$. This is
the universal commutative algebra generated by $\Z_p[G]/N_G\Z_p$, so
it is $S^*(\Z_p[G]/N_G\Z_p)$ which is isomorphic to $S^*(I)$ because
$\Z_p[G]/N_G\Z_p$ is isomorphic to $I$.  Therefore there is an exact
sequence $$0\rightarrow S^{n-1}(\Z_p[G])\rightarrow
S^n(\Z_p[G])\rightarrow S^{n}(I)\rightarrow 0.$$ Using the known
cohomology groups of the symmetric powers of $\Z_p[G]$ and the exact
sequence of cohomology groups of this exact sequence gives the
cohomology groups of $S^n(I)$ and hence the number of times that
$\Z_p$ and $I$ occur in the decomposition of $S^n(I)$.  This proves
lemma 2.7.

\proclaim Corollary 2.8.  The cohomology rings of the symmetric
algebras of indecomposable modules are give as follows. 

{\it \noindent 
$\H^*(S^*(\Z_p))$ is a polynomial algebra with one
ordinary generator of degree 1.

\noindent
$\H^*(S^*(\Z_p[G]))$  is a polynomial algebra 
with one ordinary generator of degree $p$.

\noindent
$\H^*(S^*(I))$ is the tensor product of a polynomial algebra with one
ordinary generator of degree $p$ and an exterior algebra with one
super generator of degree 1 if $p$ is odd. (If $p$ is even then
$\H^*(S^*(I))$ is a polynomial algebra on one super generator of
degree 1, but we do not need this case.) 
}

Proof. This follows from the calculation of the symmetric powers in
lemma 2.7. The product structure on the symmetric algebras can be
worked out using the formula $d(ab)=(da)b \pm a(db)$ relating the cup
product with the map $d$ from $\H^i(C)$ to $\H^{i+1}(A)$ associated to
the exact sequence $0\rightarrow A\rightarrow B\rightarrow
C\rightarrow 0$. This proves corollary 2.8.

\proclaim Lemma 2.9. 
Suppose that $A_{-1}=0$, $A_0,\ldots,A_n,A_{n+1}=0$ is a sequence of
modules and we have maps $d $ from $A_i$ to $A_{i+1}$ and $d^* $ from
$A_{i+1}$ to $A_{i}$ for all $i$ such that $d^2=0$, $d^{*2}=0$, and
$dd^*+d^*d=k$ for some fixed integer $k$ coprime to $p$. Then
$$[A_0]-[A_1]+[A_2]-\cdots +(-1)^n[A_n]=0.$$

Proof. The fact that $k$ is a unit in $\Z_p$ implies that $A_i$ is the
direct sum $dd^*(A_i)\oplus d^*d(A_i)$ and that $d$ is an isomorphism
from $d^*d(A_i) $ to $dd^*(A_{i+1})$. The lemma follows from this
because the sequence splits as the sum of these isomorphisms.  (Notice
that the relation between the $[A_i]$'s does not follow just from the
exactness of $d$; for example, look at the short exact sequence
$0\rightarrow I\rightarrow \Z_p[G]\rightarrow
\Z_p\rightarrow 0$. It is easy to check that in this case there
is a map $d^*$ with $d^{*2}=0$ and $dd^*+d^*d=p$, so the condition
that $(k,p)=1$ is also necessary.) This proves lemma 2.9. 

\proclaim Lemma 2.10. Suppose that $A$ is a finitely generated free module
over $\Z$ or $\Z_p$ with a symmetric self dual form (,). Then this
induces a symmetric self dual form (defined in the proof) on any
exterior power $\Lambda^n(A)$.

Proof. We define the bilinear form on $\Lambda^n(A)$ by
$$(a_1\wedge\cdots\wedge a_n,b_1\wedge\cdots\wedge b_n)
=\sum_{\sigma\in
S_n}\chi(\sigma)(a_1,b_{\sigma(1)})\cdots(a_n,b_{\sigma(n)})$$ where
$\chi(\sigma)$ is $1$ or $-1$ depending on whether $\sigma$ is an even
or odd permutation.  This obviously defines a symmetric bilinear form
on $\Lambda^n(A)$ with values in $\Z$ or $\Z_p$, and we have to check
it is unimodular. The determinant of $\Lambda^n(A)$ is equal to
$a\det(A)^b$ for some constants $a$ and $b$ depending on $n$, but not
on $A$. We can work out $a$ by letting $A$ be the lattice with a base
of orthogonal elements $a_1,a_2,\ldots$ with norms $\pm 1$, when we
see that $\Lambda^n(A)$ also has an orthogonal base of elements
$a_{i_1}\wedge a_{i_2}\cdots$ with $i_1<i_2<\cdots$ with norms $\pm
1$, so both $A$ and $\Lambda^n(A)$ have discriminant 1 and so $a=1$.
Hence if $A$ is self dual then $\Lambda^n(A)$ has discriminant 1 and
must also be self dual. This proves lemma 2.10.  (It is easy to check
that the constant $b$ is equal to $\dim(\Lambda^n(A))/\dim(A)$ but we
do not need this.)

The analogue of this lemma for symmetric powers is false.  The proof
breaks down because the element $a_{i_1}a_{i_2}\cdots$ with $i_1\le
i_2\le \cdots$ does not have norm $\pm 1$ if some of the $i_j$'s are
equal. In section 5 we will construct a graded algebra $h(A)$, such
that the homogeneous components of $h(A)$ have self dual symmetric
bilinear forms if $A$ does, and $h(A)\otimes \Q_p$ is a sum of tensor
products of symmetric powers of $A$.

The following corollary is purely for entertainment and is not used
elsewhere in the paper.
\proclaim Corollary 2.11. There is a 28-dimensional self dual 
positive definite lattice with no roots, acted on by
the orthogonal group  $O^+_8(\F_2)$.

Proof. By lemma 2.10  the lattice $\Lambda^2(E_8)$ is self dual, and is
obviously acted on by the automorphism group of the $E_8$ lattice
modulo $-1$, which is just $O^+_8(\F_2)$. It is easy to check that this
lattice has no roots either by calculation or by observing that there
is no suitable root lattice acted on by $O^+_8(\F_2)$. This proves
lemma 2.11.

The lattice $\Lambda^2(E_8)$ has theta function 
$$\sum_{\lambda\in
\Lambda^2(E_8)}q^{\lambda^2/2} = 1 + 2240q^{3/2} + 98280q^{4/2} +
1790208q^{5/2} + 19138560q^{6/2} +\cdots.$$ The vectors of norm 3 are
just the $2240$ vectors generating one of the 1-dimensional
sublattices of the form $\Lambda^2(L)$ where $L$ is one of the 1120
sublattices of $E_8$ isomorphic to $A_2$.

\proclaim
3.~The $\Z[1/2]$-form of the monster Lie algebra.

We find some relations between the dimensions of the Tate cohomology
groups $\H^i(g,V[1/2])$ by studying the Lie algebra cohomology of the
positive subalgebra of a $\Z_p$-form of the monster Lie algebra.
In the next section we will use these relations to compute the
cohomology groups in the case when $g$ has order at least 13.  We
recall that $V[1/2]$ is the self dual $\Z[1/2]$-form of Frenkel,
Lepowsky, and Meurman's monster vertex algebra [FLM] described in [B-R]
(but the arguments in this section work for any self dual $\Z_p$-form
of the monster vertex algebra).

This section should really be about the integral form of the monster
Lie algebra, but this has not been constructed yet because of the lack
of a construction of a self dual integral form of the monster vertex
algebra.  However, as we are only looking at the action of elements of
odd order $p$ in the monster and tensoring everything with $\Z_p$, a
$\Z[1/2]$-form is just as good as an integral form.

We briefly recall some properties of the monster Lie algebra $\m$ from
[B]. This is a $II_{1,1}$-graded Lie algebra acted on by the monster,
where $II_{1,1}$ is the 2-dimensional even self dual Lorentzian
lattice whose element $(m,n)$ has norm $-2mn$.  The piece of degree
$(m,n)$ of $\m$ is isomorphic to $V_{mn}$ if $(m,n)\ne (0,0)$ and to
$\R^2$ if $(m,n)=(0,0)$.  
The Lie algebra  $\m$ is the sum of 3 subalgebras $E$,
$F$, and $H$, where the Cartan subalgebra $H$ is the piece
of degree $(0,0)$ of $\m$, $E$ is the sum of all pieces of degree
$(m,n)$ of $\m$ for $m>0$, and $F$ is the sum of all pieces of degree
$(m,n)$ of $\m$ for $m<0$.

The $\Z[1/2]$ form on $V$ induces a self dual $\Z[1/2]$-form and hence
a $\Z_p$-form on $\m$.  From now on we will write $\m[1/2]$ for this
$\Z[1/2]$-form.  In the preprint version of this paper I implicitly
assumed that the degree $(m,n)$ piece of $\m$ was self dual and
isomorphic as a module over the monster to $V_{mn}$. At the last
moment I realized that this is not at all clear. The problem is that
although the no-ghost theorem gives an isomorphism of both spaces
tensored with $\Q$, and both spaces have $\Z[1/2]$-forms, there is no
obvious reason why this isomorphism should map one $\Z[1/2]$ form to
the other. (In fact for the fake monster Lie algebra 
the corresponding isomorphism does not preserve the integral forms.)
Fortunately in this paper we only need the following weaker statement:

\proclaim Assumption. If $m<p$ then the degree $(m,n)$ piece
of $\m\otimes \Z_p$ is self dual under
the natural bilinear form and isomorphic to $V_{mn}\otimes \Z_p$
as a $\Z_p$ module acted on by the monster. 

I believe I have a proof of this; if all goes well it will appear
in a paper provisionally titled ``The fake monster formal group''. 
Until this paper appears all the statements in sections 3 and 4 and
the modular moonshine conjectures for elements of type $pA$ $(p\ge 13)$
should have this assumption added as a hypothesis. 

I do not know if the condition $m<p$ in the assumption above is
necessary; if it is not, then this condition (usually added in parentheses)
can be missed out of some
of the following lemmas.  It seems very likely that self duality still
holds without it, but I am not so sure about the action of the monster
being the same if both $m$ and $n$ are divisible by $p$.

The monster Lie algebra has a 
Weyl vector $\rho$, i.e., a vector in
the root lattice $II_{1,1}$ such that $(\rho,\alpha) =
-(\alpha,\alpha)/2$ for every simple root $\alpha$.  This follows from
[B, theorem 7.2] which states that this simple roots are the vectors
$(1,n)$ for $n\ge 0$ or $n=-1$, with multiplicities $c(n)$. Hence the
Weyl vector is $\rho=(1,0)$.

We recall some facts about Hodge theory applied to Lie algebra
cohomology.  The self dual symmetric bilinear form on the monster Lie
algebra identifies $E$ with the dual space of $F$, and using the
Cartan involution which maps $F$ to $E$ we see that we have a
symmetric bilinear form on the Lie algebra $E$ (which is the
restriction of the contravariant form of the monster Lie algebra to
$E$) and in particular we can identify $E$ with its dual. The Lie
algebra cohomology is the cohomology of an operator $d$ which maps
each space $\Lambda^i(E)$ to $\Lambda^{i+1}(E)$.  The adjoint of this
operator is denoted by $d^*$, and we define the Laplace operator
$\Omega$ by $\Omega=dd^*+d^*d$.  As for any generalized Kac-Moody
algebra with a Weyl vector, the action of the Laplace operator
$\Omega=dd^*+d^*d$ on the degree $\alpha=(m,n)\in II_{1,1}$ piece of
$\Lambda^*(E)$ is given by multiplication by
$(\alpha,\alpha+2\rho)/2=(m-1)n$.  The Lie algebra cohomology group
$H^i(E)$ of $E$ can be identified with the zero eigenspace of $\Omega$
on $\Lambda^i(E)$. The Laplace operator restricted to $\Lambda^1(E)=E$
is half the partial Casimir operator $\Omega_0$ defined in [K, section
2.5]; the extra factor of $1/2$ is necessary to make things work well
over rings not containing $1/2$, but this is not important in this
paper as we work over the ring of $p$-adic integers for odd $p$ which
always contains $1/2$.

\proclaim Proposition 3.1. The operators $d$ and $d^*$ 
act on the $\Z[1/2]$-form $\Lambda^*(E[1/2])\subset \Lambda^*(E)$
of the exterior algebra (in all degrees $(m,n)$ with $m<p$).

Proof. The main point is that all homogeneous pieces of the
$\Z[1/2]$-form of the exterior algebra are self dual under the
contravariant bilinear form. This follows from lemma 2.10 and the fact
that each homogeneous component is a finite sum of tensor products of
exterior powers of the homogeneous components of $E[1/2]$. The result
now follows from the fact that if $d$ is a homomorphism of finite
dimensional free $\Z[1/2]$-modules from $A$ to $B$ and $A$ and $B$
have self dual symmetric inner products then the adjoint $d^*$ from
$B$ to $A$ is well defined. This proves proposition 3.1.

We let $R(0)$ be the ring of formal Laurent series $K[[r,q]][1/r,1/q]$
with coefficients in $K$, and we let $R(1)$ be the space of power
series in $R(0)$ only involving monomials $r^mq^n$ with $p|mn$, and we
let $R(2)$ be the subring of $R(0)$ of power series only involving
powers of $r^p$ and $q^p$. So we have the inclusions $R(2)\subset
R(1)\subset R(0)$ and $R(1)$ is a module over $R(2)$ (but is not a
ring itself). (We can also define $R(0)$, $R(1)$, and $R(2)$ for
arbitrary generalized Kac-Moody algebras, in particular for the fake
monster Lie algebra: $R(0)$ is a completion of the group ring of the
root lattice, $R(1)$ is a subspace corresponding to the elements
$\alpha$ of the root lattice with $p|(\alpha,\alpha)/2$, and $R(2)$ is
the subring corresponding to $p$ times the root lattice.)  If $A$ is a
$\Z[1/2]$-module then we write $[A]$ for $[A\otimes_{\Z[1/2]}\Z_p]$.

\proclaim Lemma 3.2. If $(p,(m-1)n)=1$ (and $m<p$) then the piece of degree 
$(m,n)$ of $$[\Lambda^*(E[1/2])] =
[\Lambda^0(E[1/2])]-[\Lambda^1(E[1/2])]+[\Lambda^2(E[1/2])]-\cdots$$
vanishes.  In other words $$r^{-1}\Lambda^*(\sum_{m>0,n\in \Z}r^mq^n
[V_{mn}])$$ lies in $R(1)$.

Proof. The Laplace operator $\Omega=dd^*+d^*d$ acts as multiplication
by $(m-1)n$ on the degree $(m,n)$ piece of $\Lambda^*(E[1/2])$, so the
result follows from lemma 2.9 if we let $A_i$ be the degree $(m,n)$
piece of $\Lambda^i(E[1/2])$. This proves lemma 3.2.

\proclaim Lemma 3.3. If $A$ is a $\Z_p[G]$-module with $\Tr(g|A)=c^+$,
$f([A])=c^-$, (and $m<p$) then
$$f(\Lambda^*(r^mq^nA))= (1-r^mq^n)^{(c^++c^-)/2}(1+r^mq^n)^{(c^+-c^-)/2}U$$
where $U$ is a unit of the ring $R(2)$. 

Proof. It is sufficient to prove this for an indecomposable module $A$
because both sides are multiplicative. We just check each of the three
possibilities for $A$ using lemma 2.6.

If $A=\Z_p$ then $c^+=c^-=1$ and $f(\Lambda^*(r^mq^n[\Z_p])) = 1-r^mq^n$,
so lemma 3.3 is true in this case.

If $A=\Z_p[G]$ then $c^+=c^-=0$ and $f(\Lambda^*(r^mq^n[\Z_p[G]])) =
1-r^{mp}q^{np}$ which is a unit in $R(2)$, so  lemma 3.3 is true in this case.

If $A=I$ then $c^+=-1$, $c^-=1$ and 
$$f(\Lambda^*(r^mq^n[I])) =
(1+r^{mp}q^{np})/(1+r^mq^n) = U(1+r^mq^n)^{-1},$$ so  lemma 3.3 is also
true in this case. This proves lemma 3.3.

We define integers $c^+_g(n)$ and $c^-_g(n)$ by $$\eqalign{ c^+_g(n)&=
\Tr(g|V_n)=\dim(\H^0(g,V[1/2]_n))-\dim(\H^1(g,V[1/2]_n))\cr c^-_g(n)&=
f(V[1/2]_n)=\dim(\H^0(g,V[1/2]_n))+\dim(\H^1(g,V[1/2]_n)).\cr }$$ We
use the numbers $c^+_g(n)$ and $c^-_g(n)$ rather than the apparently
simpler numbers $\dim(\H^0(g,V[1/2]_n))$ and $\dim(\H^1(g,V[1/2]_n))$
because they turn out to be coefficients of Hauptmoduls, and because
$c^+_g(n)$ is already known.

\proclaim Proposition 3.4. 
If $(p,mn)=1$ (and $m<p$) then the coefficient of $r^mq^n$ in
$$r^{-1}\prod_{m>0,n\in \Z}(1-r^mq^n)^{(c^-_g(mn)+c^+_g(mn))/2}
(1+r^mq^n)^{(c^+_g(mn)-c^-_g(mn))/2}$$ vanishes. In other words this
power series lies in $R(1)$.

Proof. We apply the homomorphism $f$ to the expression in lemma 3.2
and use lemma 3.3.  We find that the expression in this proposition is
equal to an element of $R(1)$ times a unit in $R(2)$, and is therefore
still an element of $R(1)$. Therefore its coefficients of $r^mq^n$
vanish unless $p|mn$. This proves proposition 3.4.

Proposition 3.4 can be generalized in the obvious way to any 
generalized Kac-Moody algebra which has a self dual $\Z_p$-form, 
a Weyl vector, and an integral root lattice. 

We can also ask whether or not the coefficients of $r^mq^n$
in proposition 3.4  vanish
when $p|mn$.  Some numerical calculations suggest that they usually do
not.  

\proclaim
4.~The modular moonshine conjectures for $p\ge 13$.

In this section we prove the following theorem, which completes the
proof of the modular moonshine conjectures of [R section 6] (apart
from a small technicality in the case $p=2$.)

\proclaim Theorem 4.1. If $g\in M$ is an element of prime order $p\ge 17$
or an element of type $13A$ then $\H^1(g,V[1/2])=0$. 

This implies that $^gV=\H^0(g,V[1/2])=\H^*(g,V[1/2])$ is a vertex algebra
whose modular character is given by Hauptmoduls, and whose
homogeneous components have the characters of [R definition 2]. 

The proof of this theorem will occupy the rest of this section. 
We have to show that the numbers $c^-_g(n)$ of proposition 3.4
are equal to the numbers $c^+_g(n)$, because the difference is twice the
dimension of $\H^1(g,V[1/2]_n)$. We start by summarizing 
what we know about these numbers. 
\proclaim Lemma 4.2. 
\item {1.} The numbers $c^-_g(n)$ and $c^+_g(n)$ satisfy the relations
given in proposition 3.4. 
\item {2.} The numbers $c^-_g(n)$ are integers, 
with $c^-_g(n)\equiv c^+_g(n)\bmod 2$.
\item {3.} $c^-_g(n)\ge |c^+_g(n)|$ and $(p-2)c^+_g(n)+pc^-_g(n)\le 2c^+_1(n)$,
(where $c^+_1(n)$ is the coefficient of $q^n$ in the elliptic
modular function). 
\item {4.} The numbers $c^+_g(n)$ are the coefficients of
the Hauptmodul of the element $g\in M$. 

Proof. These properties follow easily from proposition 3.4  together with the 
fact that $\Z_p\otimes V_n$ has dimension $c^+_1(n)$ and is the sum of 
$(c^-_g(n)+c^+_g(n))/2$ copies of $\Z_p$, 
$(c^-_g(n)-c^+_g(n))/2$ copies of $I$, and some copies of $\Z_p[G]$. 
This proves lemma 4.2. 

We will prove theorem 4.1 by showing that if $g$ satisfies the
conditions of theorem 4.1 then the conditions 1 to 4 in lemma 4.2 imply that
$c^-_g(n)=c^+_g(n)$.  We find a finite set of possible solutions of
conditions 1, 2, and 4 of lemma 4.2, and find that for each $p$ only
one of these satisfies the condition 3. (There are often several
solutions not satisfying condition 3.) The proof is just a long messy
calculation and the reader should not waste time looking at it.

We first put the relations into a more convenient form. 
\proclaim Proposition 4.3. There are integers $c_{m,n}$ defined for
$m,n>0$ (and $m<p$) with the following properties.
\item {1.} 
$$-\log\left(1+\sum_{m>0,n>0}c_{m,n}r^mq^n\right) \equiv
\sum_{m>0,n>0}\sum_{d|(m,n)} {c^{(-)^d}_g(mn/d^2)\over d}r^mq^n\bmod r^p$$
\item{2.} If $(p,mn)=1$ then $c_{m+1,n}=c_{m,n+1}$. 

Proof. This follows from proposition 3.4 if we multiply both sides by
$(r-q)/rq$ and then take the logarithm of both sides. This proves
proposition 4.3.

\proclaim Lemma 4.4. If $p>13$ then the values 
of $c^-_g(n)$ for $n\le 21$, $c^-_g(2n)$ for $2n\le 32$, and
$c^-_g(36)$ and $c^-_g(45)$ are given by polynomials in the numbers
$c^-_g(i)$ for $i=1,2,4,5$ with coefficients in $\Z_3$.

Proof. We can evaluate the elements $c(n), n\le 21$ and $c(2n),n\le
16$ using the argument for cases 1 and 2 of lemma 4.7 below.  We can
evaluate $c(36)$ by looking at the coefficient of $r^4q^9$ in
proposition 4.3 and using the fact that we know $c_{m,n}$ for $m\le
3$, $n\le 8$. Similarly we can evaluate $c(45)$ by looking at the
coefficient of $r^5q^9$.  Notice that the term
$(\sum_{m,n}c_{m,n}r^mq^n)^3/3$ has 3's only in the denominators of
coefficients of $r^mq^n$ when $3|m$ and $3|n$, so we do not get
problems from this for the coefficients of $r^4q^9$ and $r^5q^9$ (but
we do get problems from this if we try to work out $c^-_g(27) $ using
the same method). This proves lemma 4.4.

If $p=13$ we run into trouble when we try to determine $c^-_g(26)$;
this is why the argument in lemma 4.4 does not work in this case. Also
the argument breaks down if we try to work out $c(27)$, because when
we look at the coefficient of $r^3q^9$ we get an extra term
$c^-_g(3)/3$ which does not have coefficients in $\Z_3$. This is why
we use the coefficients $c^-_g(1)$, $c^-_g(2)$, $c^-_g(4)$, and
$c^-_g(5)$ rather than $c^-_g(1)$, $c^-_g(2)$, $c^-_g(3)$, and
$c^-_g(5)$.

Lemma 4.4 is the reason why our arguments do not work for $p<13$
(and do not work so smoothly for $p=13$). As $p$ gets smaller 
we have fewer relations to work with and it gets more difficult
to determine all coefficients in terms of the $c^-_g(m)$'s for small $m$. 
It is probably possible to extend lemma 4.4  to cover some smaller primes
with a lot of effort. Fortunately it is not necessary 
to do the cases $p<13$ because these cases have already been 
done by explicit calculations in [B-R]. 

\proclaim Lemma 4.5. If $c^-_g(1)$, $c^-_g(2)$, $c^-_g(4)$, and $c^-_g(5)$ 
are known mod $3^n$ for some $n\ge 1$, and $p> 13$, and
$(c^-_g(1),3)=1$ if $p=13$, then these values of $c^-_g(m)$ are
determined mod $3^{n+1}$.

Proof. By lemma 4.4 we see that we can determine the values mod $3^n$
of $c^-_g(n)$ for $n\le 16$, $c^-_g(2n)$ for $2n\le 32$, and
$c^-_g(36)$ and $c^-_g(45)$.

But now if we look at the coefficient of $r^3q^{3m}$ of 3.4 with
$m=1,2,4,$ or $5$, we see that $c^-_g(9m)+c^+_g(2m)/2+c^-_g(m)/3=
c^-_g(m)^3/3 + $ (some polynomial in known $c^-_g(i)$'s with
coefficients that are $3$-adic integers).  But $c^-_g(m)^3/3\bmod 3^n$
depends only on $c^-_g(m)\bmod 3^n$ if $n\ge 1$, so we can determine
$c^-_g(m)/3\bmod 3^n$ and hence $c^-_g(m)\bmod 3^{n+1}$.  This proves
lemma 4.5.

This lemma is the reason that we use 3-adic rather than 2-adic
approximation. If we try to prove the lemma above for $2^{n+1}$
instead of $3^{n+1}$, all we find is that we can determine the
$c^+_g(m)$'s mod $2^{n+1}$, which is useless because we already know
these numbers.

\proclaim Proposition 4.6. 
If the numbers $c^-_g(n)$ satisfy the conditions 1, 2, and 4 of lemma
4.2 then the numbers $c^-_g(1)$, $c^-_g(2)$, $c^-_g(4)$, and
$c^-_g(5)$ are congruent mod $3^{29}$ to the coefficients of $q$,
$q^2$, $q^4$, and $q^5$ of one of the following power series. (The
second column is a genus zero group whose Hauptmodul appears to be the
function with coefficients $c^-_g(n)$. This has not been checked
rigorously as it is not necessary for the proof of theorem 4.1.)
\vfill\line{}
$$\matrix{
p=71&\Gamma_0(71)+&
\quad&	q^{-1}	&+q	&+q^2	&+q^3	&+q^4	&+2q^5	&+\cdots\cr
p=59&\Gamma_0(59)+&
\quad&	q^{-1}	&+q	&+q^2	&+2q^3	&+2q^4	&+3q^5	&+\cdots\cr
p=47&\Gamma_0(47)+&
\quad&	q^{-1}	&+q	&+2q^2	&+3q^3	&+3q^4	&+5q^5	&+\cdots\cr
p=47&\Gamma_0(94)+&
\quad&	q^{-1}	&+q	&	&+q^3	&+q^4	&+q^5	&+\cdots\cr
p=41&\Gamma_0(41)+&
\quad&	q^{-1}	&+2q	&+2q^2	&+3q^3	&+4q^4	&+7q^5	&+\cdots\cr
p=41&\Gamma_0(82|2)+&
\quad&	q^{-1}	&	&	&+q^3	&	&+q^5	&+\cdots\cr
p=31&\Gamma_0(31)+&
\quad&	q^{-1}	&+3q	&+3q^2	&+6q^3	&+9q^4	&+13q^5	&+\cdots\cr
p=31&\Gamma_0(62)+&
\quad&	q^{-1}	&+q	&+q^2	&+2q^3	&+q^4	&+3q^5	&+\cdots\cr
p=29&\Gamma_0(29)+&
\quad&	q^{-1}	&+3q	&+4q^2	&+7q^3	&+10q^4	&+17q^5	&+\cdots\cr
p=29&\Gamma_0(58|2)+&
\quad&	q^{-1}	&+q	&	&+q^3	&	&+q^5	&+\cdots\cr
p=23&\Gamma_0(23)+&
\quad&	q^{-1}	&+4q	&+7q^2	&+13q^3	&+19q^4	&+33q^5	&+\cdots\cr
p=23&\Gamma_0(46)+23&
\quad&	q^{-1}	&	&-q^2	&+q^3	&-q^4	&+q^5	&+\cdots\cr
p=23&\Gamma_0(46)+&
\quad&	q^{-1}	&+2q	&+q^2	&+3q^3	&+3q^4	&+5q^5	&+\cdots\cr
p=19&\Gamma_0(19)+&
\quad&	q^{-1}	&+6q	&+10q^2	&+21q^3	&+36q^4	&+61q^5	&+\cdots\cr
p=19&\Gamma_0(38)+&
\quad&	q^{-1}	&+2q	&+2q^2	&+5q^3	&+4q^4	&+9q^5	&+\cdots\cr
p=19&\Gamma_0(38|2)+&
\quad&	q^{-1}	&+2q	&	&+q^3	&	&+3q^5	&+\cdots\cr
p=17&\Gamma_0(17)+&
\quad&	q^{-1}	&+7q	&+14q^2	&+29q^3	&+50q^4&+92q^5&+\cdots\cr
p=17&\Gamma_0(34)+&
\quad&	q^{-1}	&+3q	&+2q^2	&+5q^3	&+6q^4	&+12q^5	&+\cdots\cr
p=17&\Gamma_0(34|2)+&
\quad&	q^{-1}	&+q	&	&+3q^3	&	&+4q^5	&+\cdots\cr
p=13&\Gamma_0(13)+&
\quad&	q^{-1}	&+12q	&+28q^2	&+66q^3	&+132q^4&+258q^5&+\cdots\cr
p=13&\Gamma_0(26)+&
\quad&	q^{-1}	&+4q	&+4q^2	&+10q^3	&+12q^4	&+26q^5	&+\cdots\cr
p=13&\Gamma_0(26|2)+&
\quad&	q^{-1}	&+2q	&	&+4q^3	&	&+6q^5	&+\cdots\cr
}$$
It is also possible (but  unlikely) that there are other solutions for
$p=13$ for which $c^-_g(1)$ does not satisfy the inequalities
$0\le c^-_g(1)< 3^{10}$. 

Proof. For each prime $p$ we use a computer to test all 81
possibilities for $c^-_g(1)$, $c^-_g(2)$, $c^-_g(4)$, and
$c^-_g(5)\bmod 3$.  If $p>13$ then we can calculate the $p$-adic
expansion of all the coefficients $c^-_g(n)$ recursively using
lemma 4.5 above, and we reach a contradiction by looking at
coefficients of proposition 4.3 except in the cases above. For $p=13$
this does not quite work as we have not shown the values mod $3^n$
determine those mod $3^{n+1}$, so we can adopt the crude procedure of
just testing all $3^4$ possibilities for the $c^-_g(i)$'s mod
$3^{n+1}$ for each solution mod $3^n$ we have found, and checking to
see which of them leads to contradictions. There were some cases for
$p=13$ where this did not lead to a contradiction but did at least
lead to the conclusion that $c^-_g(1)<0$ or $c^-_g(1)\ge 3^{10}$. (If
$c^-_g(1)$ is not divisible by 3 then even when $p=13$ the numbers
$c^-_g(i)\bmod 3^n$ determine the numbers $c^-_g(n)\bmod 3^{n+1}$.
When $c^-_g(1)$ was divisible by 3, there were several cases where the
coefficient $c^-_g(5)$ was not determined mod $3^{n+1}$ by the
identities used by the computer program. However in all the cases
looked at by the computer with $3|c^-_g(1)$, the values of $c^-_g(1)$,
$c^-_g(2)$, $c^-_g(4)$, and $3c^-_g(5)\bmod 3^n$ uniquely determined
their values mod $3^{n+1}$.)  This proves proposition 4.6, at least if
one believes the computer calculations.  (Anyone who does not like
computer calculations in a proof is welcome to redo the calculations
by hand.)

\proclaim Lemma 4.7. If the coefficients 
 $c^-_g(1)$, $c^-_g(2)$, $c^-_g(4)$, and $c^-_g(5)$ are equal to the
coefficients $c^+_g(1)$, $c^+_g(2)$, $c^+_g(4)$, and $c^+_g(5)$ for
$g$ an element of prime order $p>13$ or an element of type $13A$ then
all the coefficients are determined by the relations of proposition
3.4.

Proof. This proof is just a long case by case check using induction.
We will repeatedly use the fact that the coefficients of $r^mq^n$ of
both sides of proposition 4.3 are equal.

Case 1. $c^-_g(2n)$ for $2n\not\equiv 0\bmod p$.  By looking at the
coefficient of $r^2q^n$ and using the fact that $c_{2,n}=c_{1,n+1}$
(as $(p,n)=1$) we see that $c^-_g(2n)$ can be written as a polynomial
in $c(n+1)$ and $c(i)$ for $1\le i <n$.

Case 2. $c^-_g(2n+1)$ for $2n+1\not\equiv \pm 1\bmod p$.  By looking
at the coefficient of $r^2q^{2n}$ we can express $c^-_g(2n+1)$ in
terms of $c^-_g(4n)$ and known quantities.  But $c^-_g(4n) $ can be
evaluated by looking at the coefficient of $r^4q^n$ and using the fact
that $c_{4,n}=c_{3,n+1}=c_{2,n+2}$ (as $2n+1\not\equiv \pm 1\bmod p$).

The recursion formulas given by cases 1 or 2 can be written out explicitly
and are given in [B, (9.1)]. In that paper they are sufficient to determine
all coefficients because the analogues of the relations in 
proposition 3.4 hold whenever $mn\ne 0$. 

The remaining cases of the proof are unusually tedious so we only sketch
them briefly.

Case 3. $c^-_g(2n)$ for $2n\equiv 0\bmod p$, $n$ odd.
We can express this in terms of $c_{2,n}=c_{3,n-1}$, which we can evaluate
by calculating $c^-_g(3(n-1))=c^-_g(2(3(n-1)/2))$ using case 1. 

For the rest of the proof it is convenient to evaluate the cluster of
elements $c^-_g(2pn-1)$, $c^-_g(2pn+1)$, and $c^-_g(4pn)$,
simultaneously. (These are the only cases left to do.)

Case 4. $2np\equiv 2\bmod 3$. We evaluate $c^-_g(2np+1)$ in terms of
$c_{3,(2np+1)/3}$. We evaluate $c^-_g(2np-1)$ by expression
$c^-_g(2np-1)c(2)$ in terms of $c_{2,2np+1}$ and evaluating the latter
by comparing it with $c_{3,2(2np+1)/3}$.  We evaluate $c^-_g(4np)$ by
expressing $c_{3,2np+2}$ in terms of $c^-_g(2)c_{2,2np}$ and
evaluating the latter by comparing it with $c_{2,3np+3}$.

Case 5. $2np\equiv 1\bmod 3$. We evaluate $c^-_g(2np-1)$ in terms of
$c_{3,(2np-1)/3}$. We evaluate $c_{1,2np+5}$ by using $c_{3,(2np+5)/3}
$ and then evaluate in turn $c_{2,2np+4}=c_{1,2np+5}$,
$c^-_g(1)c_{1,2np+3}$ (using the coefficient of $r^2q^{2np+4}$),
$c_{2,2np+3}=c_{1,2np+3}$, $c^-_g(1)c_{1,2np+1}$ which gives us the
value of $c_{1,2np+1}$. We evaluate $c_{2,2np}$ as in case 4.

Case 6. $2np\equiv 0\bmod 3$. By looking at $c_{2,2np+2}$ and
$c_{2,2np+3}$ we can find expressions for the linear combinations
$c^-_g(1)c^-_g(2np+1)+c^-_g(3)c^-_g(2np-1)$ and
$c^-_g(2)c^-_g(2np+1)+c^-_g(4)c^-_g(2np-1)$. If $p\ne 59$ or 71 then
these two linear equations are independent so we can solve for
$c^-_g(2np+1)$ and $c^-_g(2np-1)$, and then find $c_{2,2np}$ as in
case 4. If $p=59$ or 71 then an even more tedious argument using
$c_{2,2np+2}$, $c_{3,2np+2}$, and $c_{3,2np+3}$ produces 3 independent
linear relations between $c^-_g(2pn-1)$, $c^-_g(2pn+1)$, and
$c^-_g(4pn)$, so we can solve for them.

This proves lemma 4.7. 

We can now complete the proof of theorem 4.1.  If we look at the
possible solutions for the coefficients $c^-_g(n)$ given by
proposition 4.6, we see that all the solutions except those congruent
mod $3^{29}$ to the Hauptmoduls of the elements $g\in M$ are ruled out
by the inequalities of lemma 4.2. As $3^{29}>c^+_1(i)>c^-_g(i)\ge 0$
for $i\le 5$ we see that $c^+_g(i)=c^-_g(i)$ for $i\le 5$.  The
coefficients $c^+_g(n)$ also satisfy the identities of proposition 4.3
satisfied by the $c^-_g(n)$'s; for example, this follows easily from
the identity $$\sum_{n\in \Z}c^+_g(n)r^n-\sum_{n\in \Z} c^+_g(n)q^n
=r^{-1}\prod_{m>0,n\in \Z}
(1-r^{m}q^{n})^{c^+_g(mn)}(1-r^{pm}q^{pn})^{c^+_g(pmn)}$$ [B, page
434]. Therefore $c^-_g(n)=c^+_g(n)$ for all $n$ by lemma 4.7 because
they both have the same values for $n\le 5$ and both satisfy the same
recursion relations given in proposition 3.4.  This proves theorem
4.1.
 
\proclaim
5.~Calculation of some cohomology groups. 

In this section we calculate the Tate cohomology groups of some spaces
that appear in the construction of the monster vertex algebra.  In the
next section we will use these calculations to find the Tate
cohomology $\H^*(g,V[1/2])$ for elements $g$ of types $3B$, $5B$, $7B$
and $13B$.

We will write $g_0,g_1,\ldots,g_{p-1}$ for the elements $1,g,\ldots
,g^{p-1}$ of the group ring $\Z_p[G]$ to avoid confusing the group
ring multiplication with the multiplication in other rings containing
it.

Suppose that $A$ is a free $\Z_p$-module. We write $h(A)$ for the
$\Z$-graded $\Z_p$-algebra generated by elements $h_n(a)$ of degree
$n$ for $a\in A$, $n\ge 0$ with the relations
\item{1.} $ h_n(\lambda a)=\lambda^n h_n(a)$ if $\lambda\in \Z_p$. 
\item{2.} $h_0(a)=1$, $h_n(a)=0$ if $n<0$. 
\item{3.} $h_n(a+b)=\sum_{i\in \Z}h_i(a)h_{n-i}(b)$. 

If $A$ has a basis $a_1,\ldots,a_m$ then $h(A)$ is the polynomial algebra
generated by elements $h_n(a_i)$ for $n>0$, $1\le i\le m$. An action of $G$
on $A$ induces an action of $G$ on $h(A)$ in the obvious way. 

\proclaim Lemma 5.1. The cohomology ring 
$\H^*(g,h(\Z_p[G]))$ is the polynomial algebra over $\F_p$ generated
by the ordinary elements $h_n(g_0)+h_n(g_1)+\cdots+h_n(g_{p-1})$ for
$n>0$.

Proof. This follows from the fact that $h(\Z_p[G])$ is the tensor product
of the rings $A_n$, where $A_n$ is the polynomial algebra
generated by the elements $h_n(g_i)$ for $0\le i<p$. We see 
that $A_n$ is isomorphic to the symmetric algebra $S^*(\Z_p[G])$
whose cohomology is given by lemma 2.8, and we then apply corollary 2.2 to 
finish the proof. (Notice that the space spanned by the elements
$h_n(g_i)$ for $0\le i<p$ is not $h_n(\Z_p[G])$; in fact,
the latter is not even a $\Z_p$-module.)  This proves lemma 5.1.

\proclaim Proposition 5.2. 
The cohomology ring $\H^*(g,h(I))$ is isomorphic to the free exterior
algebra with one super generator in each positive degree not divisible
by $p$.  The automorphism $-1$ of $I$ multiplies each generator by
$-1$ and hence acts as $-1$ on $\H^1(g,h(I))$ and as 1 on
$\H^0(g,h(I))$.

Proof. The algebra $h(I)$ is the quotient of $h(\Z_p[G])$ by the ideal
generated by the elements $h_n(N_G)$.  We let $R_m$ be the quotient of
$h(\Z_p[G]) $ by the elements $h_n(N_G)$ for $n\le m$, and we will
prove by induction on $m$ that $\H^*(R_m)$ is the tensor product of
the free polynomial algebra generated by the ordinary elements
$h_n(g_0)+h_n(g_1)+\cdots+h_n(g_{p-1})$ (of degree $n$) for $n>m/p$
and the free exterior algebra generated by super elements of each
degree $n$ with $0<n\le m$, $(p,n)=1$ (of degree $n$).

To do this we look at the exact sequence $$0\rightarrow
R_{m-1}\rightarrow R_{m-1} \rightarrow R_{m}\rightarrow 0$$ where the
map from $R_{m-1}$ to $R_{m-1}$ is given by multiplication by
$h_m(N_G)$, and work out the cohomology of $R_m$ using the cohomology
exact sequence of this short exact sequence.  We have to consider the
cases $p|m$ and $(p,m)=1$ separately.

If $(p,m)=1$ then multiplication by $h_m(N_G)$ induces the zero map on
cohomology, so that $\H^*(g,R_m)$ is the sum of $\H^*(g,R_{m-1})$ and
$\H^*(g,R_{m-1})$ with the degree shifted by $m$ and with the ordinary
and super parts exchanged.  In other words, $\H^*(g,R_m)$ is the
tensor product of $\H^*(g,R_{m-1})$ with an exterior algebra generated
by one super element of degree $m$. This proves the induction step
when $(p,m)=1$.

If $p|m$ then multiplication by $h_m(N_G)$ induces an injective map on
cohomology, so that $\H^*(g,R_m)$ is the quotient of $\H^*(g,R_{m-1})$
by the image of $h_m(N_G)$, in other words it is just
$\H^*(g,R_{m-1})$ with a polynomial generator of degree $m$ killed
off.

Hence we see that to go from $\H^*(h(R_{m-1}))$ to $\H^*(h(R_{m}))$ we
either kill off a polynomial generator or add an exterior algebra
generator. This proves the statement about $\H^*(h(R_{m}))$ for all
$m$ by induction, and hence proves proposition 5.2.

For completeness we remark that the cohomology ring $\H^*(g,h(\Z_p))$
is a polynomial algebra over $\F_p$ with one ordinary generator of
each positive degree.

For any $\Z_p$-module $A$ we define $\omega$ to be the automorphism of
$h(A) $ taking $h_n(a) $ to $(-1)^nh_n(-a)$. We define $h_\omega(A)$
to be the largest quotient ring of $h(A)$ on which the automorphism
$\omega$ acts trivially.  (Notice that this is smaller than the
largest quotient module on which $\omega$ acts trivially.)

\proclaim Lemma 5.3. The cohomology ring $\H^*(g,h_\omega(\Z_p[G]))$
is the polynomial ring on ordinary generators of degrees $(2n+1)p$ for
$n\ge 0$.

Proof. The algebra $h_\omega(\Z_p[G])$ is the polynomial algebra
generated by the elements $h_{2n+1}(g_i)$ for $0\le i<p$, $n\ge 0$,
and lemma 5.3 follows from this.

\proclaim Lemma 5.4. The cohomology ring $\H^*(g,h_\omega(I))$
is the exterior algebra on generators of degrees $(2n+1)$ for $n\ge
0$, $(n,p)=1$. In particular it is generated (as a ring) by super
elements of odd degrees, so that $\H^0(g,h_\omega(I))$ vanishes in odd
degrees and $\H^1(g,h_\omega(I))$ vanishes in even degrees.

This follows from lemma 5.3 in the same way that lemma 5.2 follows
from 5.1.

Remark. The ring $h(\Z_p)$ can be identified with the ring of
symmetric functions over $\Z_p$ if we identify $h_n(1)$ with the
$n$'th complete symmetric function. In this case the automorphism
$\omega$ taking $h_n(1)$ to $(-1)^nh_n(-1)$ is just the automorphism
$\omega$ of [Mac, p. 14] taking the complete symmetric functions to
the elementary symmetric functions.

Remark. More generally we can form a quotient ring associated to any
automorphism $\omega$ of $A$ of order $n$ prime to $p$ such that
$h_\omega(A)$ is the quotient ring associated to the automorphism
$\omega=-1$ of $A$.  To do this we let $R_p$ be the unramified
extension of $\Z_p$ generated by a primitive $n$'th root of unity
$\zeta$. We then form the largest quotient ring of $h(A)\otimes R_p$
whose elements are fixed by the automorphism taking $h_n(a)$ to
$\zeta^nh_n(\omega(a))$. For example this space appears when we try to
construct the monster vertex algebra from a fixed point free
automorphism $\omega$ of the Leech lattice $\Lambda$.

The extraspecial group $2^{1+24}$ occurring in the centralizer of an
element of type $2B$ in $M$ has a unique $2^{12}$-dimensional
irreducible module $A$ over $\Z[1/2]$. Any element of odd order in
$\Aut(\Lambda)$ has an induced action on $2^{1+24}$ and on $A$.

\proclaim Lemma 5.5. If $g$ is an element of order $3$, $5$, 7, or 13
acting on the Leech lattice $\Lambda$ with no nonzero fixed vectors
then as a $\Z_p[G]$-module, $A\otimes\Z_p$ is the sum of $\Z_p$ and
some copies of $\Z_p[G]$.

Proof. We let $B$ be the $2^{24}$-dimensional $2^{1+24}$-module
$A\otimes A$. This is isomorphic to the group ring of $2^{24}$. The
element $g$ acts fixed point freely on the group
$2^{24}=\Lambda/2\Lambda$ because it acts fixed point freely on
$\Lambda$ and has odd order, so as a module over $G$, $B\otimes \Z_p$
is the sum of $\Z_p$ and some copies of $\Z_p[G]$, and the Tate
cohomology groups are therefore $\H^0(g,B\otimes \Z_p)=\Z/p\Z$,
$\H^1(g,B\otimes \Z_p)=0$.  As $\H^*(g,B\otimes \Z_p)=\H^*(g,A\otimes
\Z_p)\otimes \H^*(g,A\otimes \Z_p)$ we see that $\H^*(g,A\otimes
\Z_p)$ is either an ordinary or a super vector space over $\F_p$ of
dimension 1.
 
Therefore as a $G$ module, $A\otimes \Z_p$ is the sum of some copies
of $\Z_p[G]$ and either $\Z_p$ or $I$. But $\dim(A)=2^{12}\equiv 1
\bmod p$ as $(p-1)|12$ 
 which rules out the second possibility.  Therefore $\H^0(g,A\otimes
\Z_p))=\Z/p\Z$, $\H^1(g,A\otimes \Z_p)=0$. This proves lemma 5.5.

A similar argument can be used to calculate $\H^*(g,A)$ for any
element $g$ of odd prime order $p$ in $\Aut(\Lambda)$.  If $2n$ is the
dimension of the subspace of $\Lambda$ fixed by $g$ then $n$ is an
integer with $2^n\equiv \pm 1\bmod p$.  Exactly one of the two
cohomology groups $\H^*(g,A)$ vanishes and the other has dimension
$2^n$ over $\F_p$. The group $\H^1(g,A)$ vanishes if $2^n\equiv 1\bmod
p$, and the group $\H^0(g,A)$ vanishes if $2^n\equiv -1\bmod p$. There
are some cases where it is $\H^0(g,A)$ and not $\H^1(g,A)$ that
vanishes; for example, if $g$ is the element of order 5 with trace
$-1$ on $\Lambda$ then $\H^0(g,A)=0$ and $\H^1(g,A)=(\Z/5\Z)^4$.

\proclaim
6.~Elements of type $3B$, $5B$, $7B$, and $13B$.

We use the result of the previous section to calculate the cohomology
groups $\H^i(g,V[1/2])$ when $g\in M$ has type $3B$, $5B$, $7B$, or
$13B$.

\proclaim Theorem 6.1. If $g\in M$ has type $3B$, $5B$, $7B$, 
or $13B$ and $\sigma $ is the element of order 2 generating the center
of $C_M(g)/O_p(C_M(g))$ then $O_p(C_M(g))$ acts trivially on $^gV$ and
$\sigma$ fixes the ordinary part of $^gV$ and acts as $-1$ on the
super part.

Proof.  We prove this by calculating the action of $\sigma$ on
$\H^*(g,V[1/2])$ explicitly. We recall from [B-R] that the
$\Z[1/2]$-form $V[1/2]$ of Frenkel, Lepowsky and Meurman's monster
vertex algebra [FLM] is the sum of two pieces $V^0$ and $V^1$ which
are the $+1$ and $-1$ eigenspaces of an element $\sigma$ of type $2B$.
(Unfortunately the details for the construction of this $\Z[1/2]$-form
have not been published anywhere.)  We now show that
$\H^0(g,V^1)=\H^1(g,V^0)=0$.

The space $V^0$ is the subspace of $V_\Lambda[1/2]$ fixed by the
automorphism induced by $-1$ acting on the Leech lattice $\Lambda$,
where $V_\Lambda[1/2] = V_\Lambda\otimes \Z[1/2]$, and $V_\Lambda$ is
the integral form of the vertex algebra of the Leech lattice,
described in [B-R].  As in [B-R], the space $V_\Lambda$ is the tensor
product of a twisted group ring $\Z[\hat \Lambda]$ and a space
isomorphic to $h(\Lambda)$. The cohomology $\H^*(g,\Z[\hat\Lambda])$
is just $\Z/p\Z$ as in [B-R] because $g$ acts fixed point freely on
$\Lambda$.  By proposition 5.2, the element $-1$ of
$\Aut(\hat\Lambda)$ acts as $+1$ on $\H^0(g,h(\Lambda))$ and as $-1$
on $\H^1(g,h(\Lambda))$. Therefore $\H^1(g, V^0)$ vanishes because it
is isomorphic to the subspace of $\H^1(g,h(\Lambda))$ fixed by $-1\in
\Aut(\Lambda)$.

The space $V^1$ is the subspace of elements of integral weight of the
tensor product $A\otimes h_\omega(\Lambda\otimes \Z_p)$, where $A$ is
the $2^{12}$ dimensional module for $2^{1+24}$ and
$h_\omega(\Lambda\otimes \Z_p)$ has its grading halved and then shifted
by $3/2$.  By lemma 5.5, $\H^*(g,A)=\F_p$, so $\H^*(g,A\otimes
h_\omega(\Lambda\otimes \Z_p))=
\H^*(g,h_\omega(\Lambda\otimes \Z_p))$. By lemma 5.4 the cohomology ring 
$\H^*(g, h_\omega(\Lambda\otimes \Z_p))$ 
is generated by super elements of
odd degree, so $\H^1(g, h_\omega(\Lambda\otimes \Z_p))$
vanishes in even degrees
and $\H^1(g, h_\omega(\Lambda\otimes \Z_p))$ vanishes in odd degrees. 
Therefore $\H^0(g,A\otimes h_\omega(\Lambda\otimes \Z_p))$ 
vanishes in integral degrees, so $\H^0(g,V^1)=0$. 

This shows that the element $\sigma$ acts as 1 on
$\H^0(g,V[1/2])=\H^0(g,V_\Lambda)=\H^0(g,h(\Lambda\otimes \Z_p))$ and
as $-1$ on $\H^1(g,V[1/2])=\H^1(g,A\otimes h_\omega(\Lambda\otimes
\Z_p))$.

Finally we show that $O_p(C_M(g))$ acts trivially on $^gV$.  The
element $\sigma$ acts trivially on the center $\langle g\rangle$ of
the extraspecial group $p^{1+24/(p-1)}$ and acts as $-1$ on the
quotient $p^{24/(p-1)}$. The element $g$ acts trivially on $^gV$ so we
get an action of $p^{24/(p-1)} $ on $^gV$.  As $\sigma$ acts as 1 on
the even part of $^gV$ and as $-1$ on the odd part, it commutes with
all automorphisms of $^gV$, so the group $p^{24/(p-1)}$ must act
trivially on $^gV$.  This proves theorem 6.1.

The group $C_M(g)$ has structure $3^{1+12}.2.Suz$, $5^{1+6}.2.HJ$,
$7^{1+4}.2.A_7$, or $13^{1+2}.2.A_4$, so we get modular vertex
superalgebras acted on faithfully by the groups $2.Suz$, $2.HJ$,
$2.A_7$, and $2.A_4$ whose modular characters are given by
Hauptmoduls.  We can calculate the characters of both the ordinary and
super parts of these superalgebras by using the fact that the
characters of their sum and difference are both given in terms of
known Hauptmoduls.

It seems likely that if $g$ has type $2B$ then the vertex superalgebra
acted on by $2^{24}.Co_1$ constructed in [B-R] is acted on trivially by
the $2^{24}$ subgroup, but this has not yet been proved. The proof
given above for the corresponding fact for elements of odd prime order
fails for order 2 because an automorphism mapping each element of
$2^{24}$ to its inverse is still the identity automorphism.

As an  example, we work out 
the dimensions of the ordinary and super parts of $^gV$ 
in the case when $g$ has type $3B$. 
The function $\sum_n\Tr(1|^gV)$ is equal to 
the Hauptmodul of the element $g$ of type $3B$ in the monster,
which is
$$q^{-1}+54q-76q^2-243 q^3+1188q^4-1384q^5-\cdots,$$
and the function $\sum_n\Tr(\sigma|^gV)$ is equal to 
the Hauptmodul of the element $\sigma g$ of type $6B$ in the monster,
which is
$$q^{-1}+78q+364q^2+1365q^3+4380q^4+12520q^5+\cdots$$
so we find that the dimension of the ordinary part of $^gV_n$
is the coefficient of $q^{n}$ of
$$q^{-1}+66q+144q^2+561q^3+2784q^4+5568q^5+\cdots$$
and  the dimension of the super part of $^gV_n$
is the coefficient of $q^{n}$ of
$$12q+220q^2+804q^3+1596q^4+6952q^5\cdots.$$
These series can be written as sums of dimensions
of ordinary representations of the groups $3.Suz$ and $6.Suz$ as
$$q^{-1}+66q+(1+143)q^2+(429+66+66)q^3+(66+66+78+429+2145)q^4+
(4\times 1+4\times 143+2\times 780+3432)q^5+\cdots$$
and
$$12q+220q^2+(780+12+12)q^3+(780+780+12+12+12)q^4+
(4\times 220+2\times 572+4928)q^5+\cdots.$$
(Some of the decompositions are ambiguous; for example $1+143$
could be $66+78$. We have chosen $1+143$ instead of $66+78$
for a reason that will appear in a moment.)
At first sight this seems to rule out a vertex algebra 
structure on any lift to characteristic 0, as the 
homogeneous components $V_n$, $V_{n+1},\ldots$ have to increase. 
However we can get around this by rescaling the grading
by a factor of $1/3$. The data above seem consistent
with the existence of a ${1\over 3}\Z$-graded vertex superalgebra
$W=\sum_{3n\in \Z}W_n$ such that the mod 3 reduction
of $W_n$ is $V_{3n}$. The vertex superalgebra $W$ is probably
acted on by $6.Suz$, with the element of order 3 in the center
acting as $\omega^{3n}$ on $W_n$. 

There is a second modular vertex super algebra acted on by each of the
groups $2.Suz$, $2.HJ$, $2.A_7$, and $2.A_4$, which can be constructed
as $\H^*(g,V_\Lambda)$ where $g$ is a fixed point free automorphism of
order 3, 5, 7, or 13 of the Leech lattice $\Lambda$. The ordinary part
of this superalgebra is the same as the ordinary part of the
superalgebra constructed from the monster vertex algebra, but its
super part is different; for example, its degree 0 piece does not
vanish.

\proclaim
7.~Open problems and conjectures.

For more open problems in this area see section 7 of [B-R].

At the end of section 6 we saw that the vertex superalgebra of an
element of type $3B$ looks as if it might be the reduction mod 3 of a
superalgebra in characteristic 0. More generally something similar
seems to happen for some other elements of the monster, as described
in the following conjecture.

\proclaim Conjecture. 
For each element $g \in M$ of order $n$ there is a
${1\over n}\Z$-graded  super module $^g\hat V=\sum_{i\in{1\over n}\Z}V_i$
over the ring $\Z[e^{2\pi i/n}]$
with the following properties.
\item{1.} $^g\hat V$ is often a vertex superalgebra
with a conformal vector and
 a self dual invariant symmetric bilinear form.
(For example, it seems to be a vertex algebra when $g$
is of type $1A$, $2B$, $3B$, $5B$, $7B$, $13B$, or $3C$,
 but possibly not if $g$ has type
$2A$ or $3A$.)
\item{2.} $^1\hat V$ is  a self dual integral form of the monster
vertex algebra (probably the one conjecturally described in 
lemma 7.1 of [B-R]).
\item {3.} If $g$ has prime order then the reduction
mod $p$ of $^g\hat V$ is the modular vertex superalgebra
$^gV$ associated to $g$ (at least when its grading is multiplied
by $n$.)
\item{4.} If all coefficients of the Hauptmodul
of $g$ are nonnegative then the super part of $^g\hat V$ vanishes.
 In this case $^g\hat V $
should be the same as the ``twisted sector''  of $g$ that appears in 
Norton's generalized moonshine conjectures [N]. 
(But notice that if the Hauptmodul of
$g$ has negative coefficients then $^g\hat V$ is certainly not the 
twisted sector of $g$; for example, $^g\hat V$ is a super vector space and 
has a nontrivial piece of degree 0.)
\item{5.} $^g\hat V$ is acted on by a central extension $(\Z/n\Z).C_M(g)$
of
$C_M(g) $ by $\Z/n\Z$. This central extension 
contains an element of order 1 or 2 acting as $-1$ on the super part
of $^g\hat V$ and as $1$ on the ordinary part. 
\item{6.} If 
$g$ and $h$ are commuting elements of $M$ of coprime orders
and $\hat g$ 
is the   lift of $g$  to the central
extension of $C_M(h)$ with the same order as $g$, 
then 
$$\Tr(\hat g|^{h}\hat V)=\Tr({gh}|^1\hat V).$$
(This identity is often
 false if the orders of $g$ and $h$ are not coprime.)
\item{7.} If $g$ and $h$ commute then
$\H^*(\hat g,^h\hat V)$ is (essentially)
the reduction mod $|g|$ of $^{gh}V$. (Actually this
is not quite correct as stated, because the reduction mod $|g|$
will be defined over some nilpotent extension of $\Z/|g|\Z$,
and we should then take the tensor product over this ring with $\Z/|g|\Z$.)

These conjectures are similar to Norton's generalized moonshine
conjectures [N] because both associate some graded space with each
element $g$ of the monster.  However they treat the problem of
negative coefficients in the Hauptmodul $T_g(\tau)$ in a different
way: Norton's conjectures use the fact that $T_g(-1/\tau)$ has
nonnegative coefficients, which are supposed to be the dimensions of
some vector spaces, while the conjectures above allow super vector
spaces which may have negative super dimension.  Perhaps both Norton's
conjectures and the conjectures above are special cases of some
conjecture which associates a $\Z/n\Z\times \Z$-graded super space
acted on by a central extension $\Z/n\Z.C_M(g)$ to each element $g\in
M$ of order $n$.

We finish by listing a few other questions.

\item{1.} Is the group $C_M(g)/O_p(C_M(g))$ the full automorphism group
of the modular vertex algebra $^gV$? If $g$ has type $3C$ it seems
possible that it is not, and the full automorphism group may be
$E_8(\F_3)$. 
Notice that if $g$ has type $2B$ then it is not yet known
that $O_2(C_M(g))$ acts trivially on $^gV$.

\item{2.} Replace the revolting argument of section 4 that involves computer
calculations with a more conceptual proof.  More generally, suppose we
have coefficients $c(n)$ of some function such that the coefficient of
$r^mq^n$ in $$r^{-1}\prod_{m>0,n\in \Z}(1-r^mq^n)^{c(mn)}$$ vanishes
whenever $(mn,N)=1$ for some fixed integer $N$. What does this imply
about the coefficients $c(n)$? (There are many solutions to this given
by taking the $c(n)$'s to be the coefficients of some Hauptmodul of
level $N$.)

\item{3.} What happens when we apply the
argument of section 3 to some automorphism of prime order of the fake
monster Lie algebra?  Do the infinite products we get represent
interesting functions?

\item{4.} In [B section 10] it was suggested that there 
should be Lie superalgebras associated with elements $g$ of the
monster, but it was not clear what the dimensions of homogeneous
components of these superalgebras should be in the cases when the
coefficients of the Hauptmodul of $g$ were neither all positive nor
alternating in sign. The results of section 6 suggest what to do in
this case: the sum of the dimensions of the ordinary and super parts
should also be given by the coefficients of another Hauptmodul.

\proclaim References.

\item{[BBCO]} C. Batut, D. Bernardi, H. Cohen, M. Olivier,
``User's guide to PARI-GP''. This guide and the PARI programs 
can be obtained by anonymous ftp from hensel.mathp7.jussieu.fr
in the directory /dist/pari.
\item{[B]}{R. E. Borcherds, 
Monstrous moonshine and monstrous Lie superalgebras, Invent. Math.
109, 405-444 (1992).}
\item{[B-R]}{R. E. Borcherds, A. J. E. Ryba, 
Modular Moonshine II, Duke Math Journal Vol. 83 No. 2, 435-459, 1996.  }
\item{[CG]} C. J. Cummins, T. Gannon, Modular equations and the genus
zero property of moonshine functions, Invent. Math. 129, 413-443 1997. 
\item{[DM]}C.~Dong and G.~Mason,
On the construction of the moonshine module as a
$\Z_p$-orbifold,
Mathematical aspects of conformal and topological field theories
and quantum groups (South Hadley,
MA, 1992), 37--52, Contemp. Math., 175, 
Amer. Math. Soc., Providence, RI, 1994.
\item{[FLM]}{I. B. Frenkel, J. Lepowsky, A. Meurman, 
Vertex operator algebras and the monster, Academic press 1988. }
\item{[Mac]} I. G. Macdonald, ``Symmetric functions and Hall polynomials'', 
Oxford University press, 1979. 
\item{[Ma]} Y. Martin, On modular invariance of completely replicable 
functions. Moonshine, the Monster,
and related topics (South Hadley, MA, 1994), 263--286, 
Contemp. Math., 193, Amer. Math. Soc., Providence, RI,
1996. 
\item{[M]} P. Montague, Third and Higher Order NFPA Twisted
Constructions of Conformal Field Theories from Lattices, 
Nuclear Phys. B 441 (1995), no. 1-2, 337--382.
\item{[N]}{S. P. Norton, Generalized Moonshine,
  Proc. Symp. Pure Math. 47 (1987) p. 208-209.}
\item{[R]}{A. J. E. Ryba, Modular Moonshine?, 
In ``Moonshine, the Monster, and related topics'', 
edited by Chongying Dong and Geoffrey Mason. 
Contemporary Mathematics, 193. American Mathematical Society,
Providence, RI, 1996. 307-336. }
\bye